# 6-Fold Quasiperiodic Tilings With Two Diamond Shapes

Theo P. Schaad

Abstract:

A set of tiles for covering a surface is composed of two types of tiles. The base shape of each one of them is a diamond or rhombus, both with angles 60 and 120 degrees. They are distinguished by labeling one as an acute diamond with a base angle of 60 degrees, the other one as an obtuse diamond with a base angle of 120 degrees. The two types of tiles can be marked with arrows, notches, or colored lines to keep them distinct. Notches can be used as matching rules such that the acute diamonds can form a star with 6-fold rotational symmetry among other matches. Similarly, three obtuse diamonds can be matched with 3-fold rotational symmetry to form a hexagon among other possibilities. Two variations of an aperiodic inflation scheme are disclosed to match nine tiles into two larger tiles. These two larger tiles of the second generation are the new base shapes following the same matching rules as the original tiles. The inflation can thus be repeated indefinitely creating an arbitrarily large surface covered with a 6-fold quasiperiodic tiling consisting of only two diamond shapes. The notches of the tiles can be creatively deformed to make Escher-esque figures.

Background:

This work is a continuation of the discovery of Roger Penrose in 1974 of a 5-fold aperiodic tiling consisting of only two distinct tiles[5]. His preferred base shapes were kites and darts, but he found that other shapes related to the kites and darts could also work, among them an unusually simple pair: two rhombs with acute angles of 72 and 36 degrees, respectively. Five-fold symmetry is not an allowed periodic tessellation, yet Penrose's aperiodic tiling has plenty of five-fold stars and other five-fold structures. The term quasiperiodic is sometimes used in this context. From a physical view, the term quasicrystal describes matter following the matching rules of quasiperiodic patterns. Another intriguing element of the Penrose pattern is that the areas of the two tiles are in the golden ratio, an irrational number that shows up in some of the most aesthetic and pleasing contexts of architecture, geometry, and nature. In the 1977 Scientific American article[4], M. Gardner asked the poignant question: Are there pairs of tiles not based on the golden ratio that force nonperiodic tiling?

R. Penrose's work was preceded by sets of square tiles, now called Wang[2] tiles, which are technically rhombs. R. Berger[3] discovered in 1966 the first aperiodic set consisting of 20,426 Wang tiles, later reduced to much smaller sets.

An 8-fold aperiodic tiling is now attributed to Ammann-Beenker[6,8]. It consists of two generalized rhombs, one a square, the other a diamond shape with an acute angle of 45 degrees. This 8-fold tiling is closely associated with the scaling factor (1 + √2), again an irrational number, sometimes referred to as the silver ratio. Since there are no 8-fold periodic tessellations, the Ammann-Beenker tiling is also aperiodic consisting of only two rhombs. The base shapes, however, do not force an aperiodic pattern without matching rules, as the square shape can be tiled in a periodic pattern.

Furthermore, 10-fold quasicrystal tilings were discovered by R. Ingalls[9] and others, consisting of only two rhombs, the same rhombs as originally disclosed by R. Penrose but with a different tiling that includes ten-fold stars. Again, there is no 10-fold periodic tessellation, so the Ingalls tiling is aperiodic. A recent paper by Gregory R. Maloney[11] describes a method for constructing planar rhombic tilings, including a novel 11-fold tessellation, generally using more than two rhomb shapes. It also contains a more extensive reference list on the subject.

There are many ways to create aperiodic tilings from periodic tessellations. For instance, one could color a periodic tile set randomly with two colors. However, finding matching rules with only two types of rhombs belongs to a special class of quasiperiodic tile sets. Such a set is described below.

## Matching rules of the 6-fold quasiperiodic tile set:

There is a periodic tessellation that consists of diamond shapes with angles of 60 and 120 degrees. It forms a subdivision of the hexagonal tessellation shown in Fig. 2 and it is also known as the rhombille tiling of the Euclidian space.

To construct a 6-fold quasiperiodic tiling with two distinct types of diamonds, they can be distinguished by color. In order to fit together with matching rules, they should have equal areas and matching edges, and should be superimposed on the hexagonal lattice without being periodic. One of the tiles coined the acute tile should have a base angle of 60 degrees and could be tiled into a 6-fold star with 6-fold rotational symmetry. Such a shape is shown in Fig. 1. Either acute angle distinguished by an arrow could be matched into a 6-fold star, but either side could be unique such that they cannot match the other side. Two embodiments of matching rules are shown in Fig. 1, one with colored bars, the other with notches. The arrow is optional since one could simply remember that the direction of the acute diamond is green-to-red bars.

It is shown in Fig. 2, that to continue tiling the stars, obtuse diamonds are needed with colored bars or notches as shown in Fig. 1. Thus, a basic set is created consisting of two types of rhombs with notches and bars that are similar to the Penrose rhombs but with 6-fold stars instead of 5-fold stars.

To show that the a tiling can cover a surface indefinitely, Penrose proposed a scheme of inflation or deflation. Such an inflation scheme is shown in Fig. 3 for both embodiments of the tiles of Fig. 1. As can be seen, all the colored bars or notches match. The new shapes resemble the previous shapes and are now the new distinct tiles of a new generation, or generation 2. Each new tile consists of nine tiles of the previous generation. The new inflated acute tile consists of 5 smaller acute and 4 obtuse tiles. Likewise, the inflated obtuse tile consists of 4 smaller acute and 5 obtuse tiles. In any larger tessellation, the number of distinct diamonds may not be exactly the same, but for a large set, the ratio will approach unity. In this set, the golden or silver ratio of the 5-fold and 8-fold rhomb areas is a simple integer (unity). The area of the diamond with unit length is an irrational number ($\sqrt{3}/2$), but the inflation scale is the integer 3, not an irrational number.

The inflation (or substitution) scheme includes obtuse diamonds that are cut in half and must be matched with the other half as the tiling continues. This is assured since it was postulated that the tiles are superimposed on a periodic rhombille tessellation. The overlapping obtuse diamonds therefore become periodic elements in the next generation. Yet the overall scheme is aperiodic and cannot be derived by translational symmetry. Higher generations are shown in Fig. 4. It can be seen that the inflation scheme holds and that a surface could be tiled in a quasiperiodic fashion indefinitely.

In what follows are different variations and embodiments of the same 6-fold scheme. A somewhat fanciful version is shown in Fig. 5 with deformed diamond shapes akin to M.C. Escher's periodic animals. In these figures, a "Fleur-de-lis" is actually a deformed obtuse diamond, but is shown in the vertical direction. Likewise the heart ("Coeur d'Alene") motif is a deformed acute diamond but shown horizontally to conform with our normal perception of these shapes. In this example, a mirror symmetry of the basic tile was retained, such that there is only one type of notches, or, equivalently, one color bar only.

To get a feel of what higher generations look like, the tiling is shown in Fig. 6 & 7 with colored gold and purple diamonds without the notches, arrows, or colored bars, but with the same inflation scheme.

Variation 2:

In the first inflation (or substitution) scheme of Fig. 3, the center tile in the next generation has an arrow with opposite (down) polarity than the new implied arrow of generation 2 (up). This appears to be somewhat arbitrary. In fact, the direction can be reversed as shown in Fig. 8 leading to a slightly different inflation scheme and a slightly different 6-fold quasiperiodic tiling as shown in Fig. 9 (compare with Fig. 4). To find two variations in nature may not be totally unexpected. For example, it is encountered with right and left-handed chirality, or, in another case, twinning quartz crystals. It would be a challenge to find out how nature would choose one of the variations, or how the researcher could distinguish the two in a physical sample.

## Aperiodicity of the 6-fold tiling:

At closer inspection of larger tilings, for instance the one in Fig. 7, clusters of tiles are noticed with perfect 6-fold rotational symmetry. Some of these clusters tend to repeat in a hexagonal lattice. For instance, at the lowest scale, there are 6-fold stars formed with acute base diamonds that repeat in a hexagonal pattern. In Fig. 7, there is a very long ribbon circling the outer parts of the tiling with striking 6-fold rotational symmetry. At the same time, in each generation, there are features that span the length scale of the inflated diamond and are clearly aperiodic relative to the underlying hexagonal lattice. A similar phenomenon appears in the 5-fold and 8-fold aperiodic tilings. There, too, 5-fold and 8-fold stars appear in the next generation followed by larger 5-fold and 8-fold structures appropriate to the inflated length scale of the higher generations.

To illustrate the appearance of these 6-fold structures, five colors were chosen in Fig. 10. All obtuse diamonds are yellow. In generation 2, the acute diamonds at the 60 degree corners are pink. These turn into pink stars in the higher generations that are all superimposed on the underlying hexagonal lattice of the second generation. In generation 3, the acute diamonds near the 60 degree corners are colored green. These turn into 6-fold rings at the higher generation and are superimposed on a hexagonal lattice of the third generation. One could easily continue this process with each generation. For instance, the short blue ribbons near the 60 degree corners in generation 4 of Fig. 10 could be assigned a new color leading to a larger six-fold structure in the next generation that repeats in a hexagonal pattern.

Returning to generation 3 in Fig. 10, the long ribbon consisting of acute base diamonds is painted blue. The remaining acute diamonds are painted red leading to aperiodic red stars in the next generation. In generation 4, four of these red stars are near the 60 degree corners and will be incorporated into a 6-fold structure in the next generation. An ever longer blue ribbon appears in generation 4 and is clearly aperiodic. The two long blue ribbons in the acute and obtuse types of inflated diamonds neatly match each other like one of the colored bars in Fig. 1. It can readily be imagined that if a huge floor were tiled with this 6-fold pattern, there will be a similar aperiodic blue ribbon of acute diamonds at the length scale of the highest generation, i.e. the size of the room, with smaller 6-fold periodic features at the appropriate length scales of the lower generations, down to periodic pink stars at the scale of the second generation.

## Darts, kites, and shields:

R. Penrose had originally found that the rhombic shapes of the 5-fold quasiperiodic tiling could be mapped into shapes that he labeled darts and kites. Conversely, any pattern with darts and kites can be transformed back into a tiling consisting of rhombs. It was therefore interesting to study if the same holds true for the 6-fold pattern. A dart can be constructed with two half-tiles (split along the long diagonal of the acute diamond). A kite can be formed with a split acute tile and an obtuse diamond. The rhombs in a tiling of the next generations can then be substituted with darts and tiles. However, that leaves untileable areas that look like "shields", or hexagons. These shields consist of a three obtuse diamonds. Curiously, the two basic diamonds form three new figures that become a set of three basic tiles that can tile a large surface with a quasiperiodic tiling.

The set of darts, kites, and shields is shown in Fig. 11. Below, the substitution rule is shown that leads to the next generation of the base shapes, a process that can be repeated indefinitely. Again, there is another variation possible that is not shown. In Fig. 12, the set of darts, kites, and shields is used to tile the basic set of diamonds. These diamonds can then be inflated as before. Half-darts, half-kites, and half-shields at the border must be matched to adjacent tiles to complete them.

## Hexagrid dual:

The multigrid method was originally developed by de Bruijn[10]. The grid consists of parallel lines that are normal to the rhomb edges. If the lines are properly joined across the tile edges, continuous straight lines are formed called Ammann lines. These lines also decorate the tiles and provide another way to distinguish the acute and obtuse diamond shapes and show the matching rules. As the grid consists of infinite parallel lines, they determine the short-range and long-range order of the tiling. If in a multigrid at most two lines intersect at the same point, it is possible to reconstruct a tiling consisting of rhombs with edges normal to the grid lines. This unique mapping makes the multigrid a dual of the tiling. The multigrid pattern is itself an aperiodic pattern if the tiling is aperiodic.

The principal directions in a hexagonal pattern are 6-fold, but the normal directions are only 3-fold. It is therefore expected that the hexagrid dual has parallel lines at 0, 120, and 240 degrees only. It was observed that if the hexagrid consists of equally spaced parallel lines, the dual is the periodic rhombille hexagonal tiling. It was also noted that a hexagrid with two alternate spacings AB is the dual of a tiling that consists of only two tiles, but is periodic (a hexagonal pattern with 6-fold stars and hexagons).

A hexagrid dual of the quasiperiodic 6-fold tiling is shown in Fig. 13. It consists of three different spacings A,B, and C. From the drawing A=1+2δ, B=2-2δ, C=3/2 in units of the rhomb edge length. δ is the distance to the nearest vertex, but it can be an arbitrary number between 0 and ½ (it simply moves the line decoration of the tiles). If δ=1/2, all spacings are identical and the rhombille periodic tiling is mapped. If δ is changed to (1/2-δ), the A and B spacings are switched. This can be seen as a switch in the polarity sign of the tile without changing the overall matching rules.

The hexagrid spacing sheds some light why there are two different variations. In variation 1 the grid spans the inflated acute tile as X=aBCa where a=A/2. If the direction of the acute tile is reversed, the spacing is Y=aCBa. Substituting three tiles along the long diagonal of the acute tile requires the sequence X → XYY. As the hexagrid is reversed on the opposite side of the 6-fold center (marked with a circle), Y → XXY. This is the essence of the substitution rule of variation 1. In the next generation we get XYY → XYYXXYXXY. The polarity of the tiles can also be changed by interchanging A with B, i.e., Y=bACb where b=B/2. This change does not change the substitution rule of variation 1. In the second drawing of Fig. 13, the substitution is Y → YXX and the mirror reverse is X → YYX. This pattern is a dual of the first variation.

|  | Substitution rule | Mirror | Dual | Dual mirror |
| --- | --- | --- | --- | --- |
| Variation 1 | X → XYY | Y → XXY | Y → YXX | X → YYX |
| Variation 2 | Y → YYX | X → YXX | X → XXY | Y → XYY |
| Periodic | X → XXX | Y → YYY | Y → XXX | X → YYY |
| Periodic | X → XYX | Y → YXY | X → YXY | Y → XYX |

As summarized in the table, there are 8 solutions (including mirror reflections) to fill an acute tile with 3 smaller ones along the long diagonal. Of these solutions, there are fundamentally two versions with their duals that form the two variations proposed with two diamond shapes. The scale factor 3 is the lowest number that leads to the possibility of an aperiodic substitution. It should be possible to extend the scheme into higher (integer) scale factors, thus this class of 6-fold aperiodic tilings could be infinitely large.

## Summary:

Two variations of a 6-fold quasiperiodic tiling were disclosed. The base tiles consist of two distinct diamonds, both with 60 and 120 degree corners. They are distinguished as an acute and obtuse diamond corresponding to the corner that can be matched into 6-fold stars or 3-fold hexagons. Matching rules are given with colored bars or notches. An inflation or substitution scheme is proposed that combines nine base shapes into larger versions of an acute and obtuse diamond shape. The process of inflation can be repeated until an arbitrarily large surface is covered in a quasiperiodic tiling. The notches of the matching rules can be creatively changed to make Escher-like figures.

This 6-fold quasiperiodic tiling is in a special class of infinite aperiodic tilings with only two distinct rhombs as base tiles. In many ways it is strikingly similar to the 5-fold Penrose tiling, in particular with the matching rules. In both patterns, the matching rules involve 2-colored bars or two type of notches. However, there is a caveat. The rhombs of the 5-fold tiling only tile in an aperiodic fashion. The obtuse diamond of the 6-fold could be tiled into a periodic hexagonal tiling, and the 6-fold stars of the acute base diamonds could be turned into hexagons with the obtuse diamonds and then tiled into a periodic hexagonal tiling. In that sense, nature may have a choice to form a periodic hexagonal or a 6-fold quasiperiodic pattern with the same base tiles.

Hermann Weyl[1] put forward the notion that symmetry equates to harmony of proportions. Irregardless if nature were to imitate art and a physical 6-fold quasicrystal were found, the quasiperiodic pattern itself adds another level to the symmetric harmony of the hexagonal tiling. The result is unique and intriguing.

The author would like to acknowledge helpful discussions with Robert Ingalls who shares a fascination with tessellations.

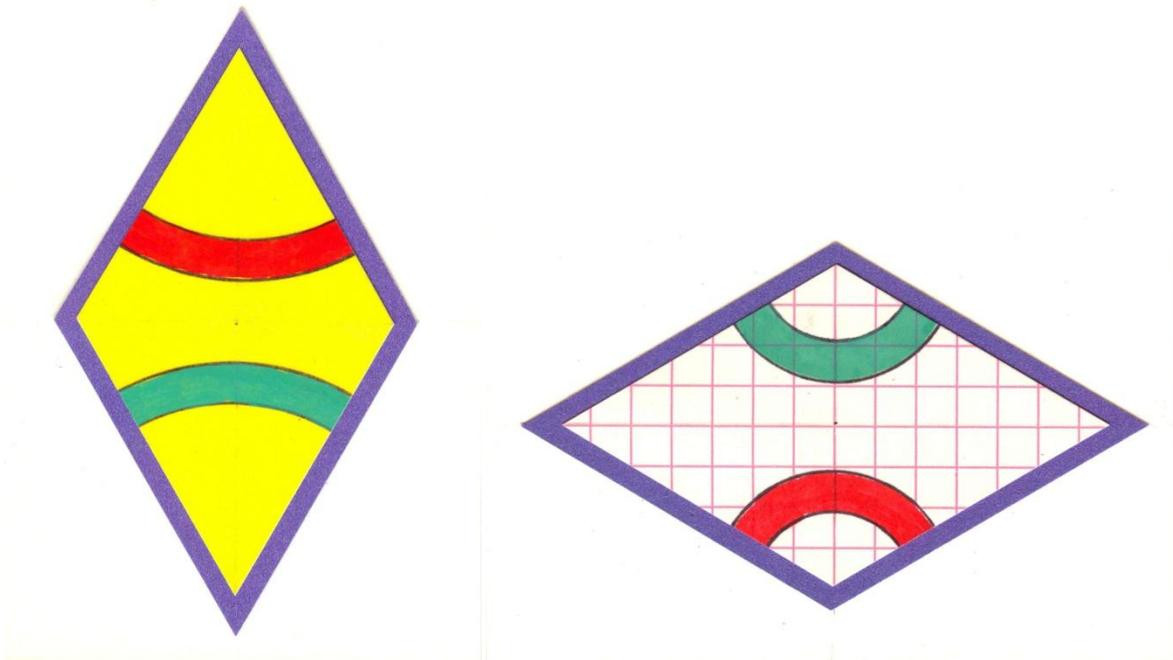

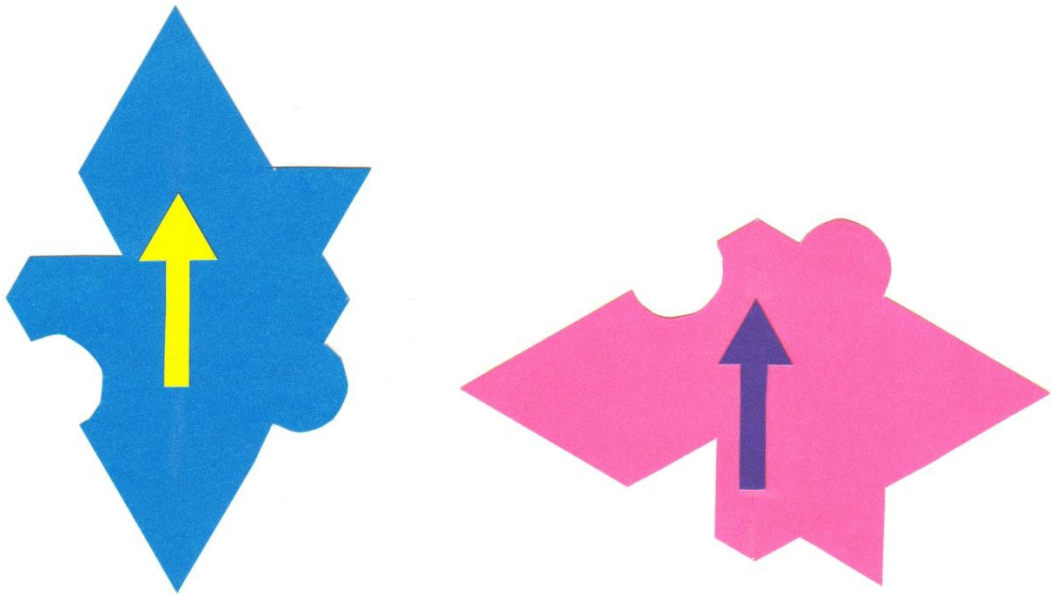

Fig. 1: Two embodiments of the acute and obtuse diamond base shapes; on top with 2-colored bars spanning the acute and obtuse angles; at bottom with notches and arrows to distinguish the two shapes.

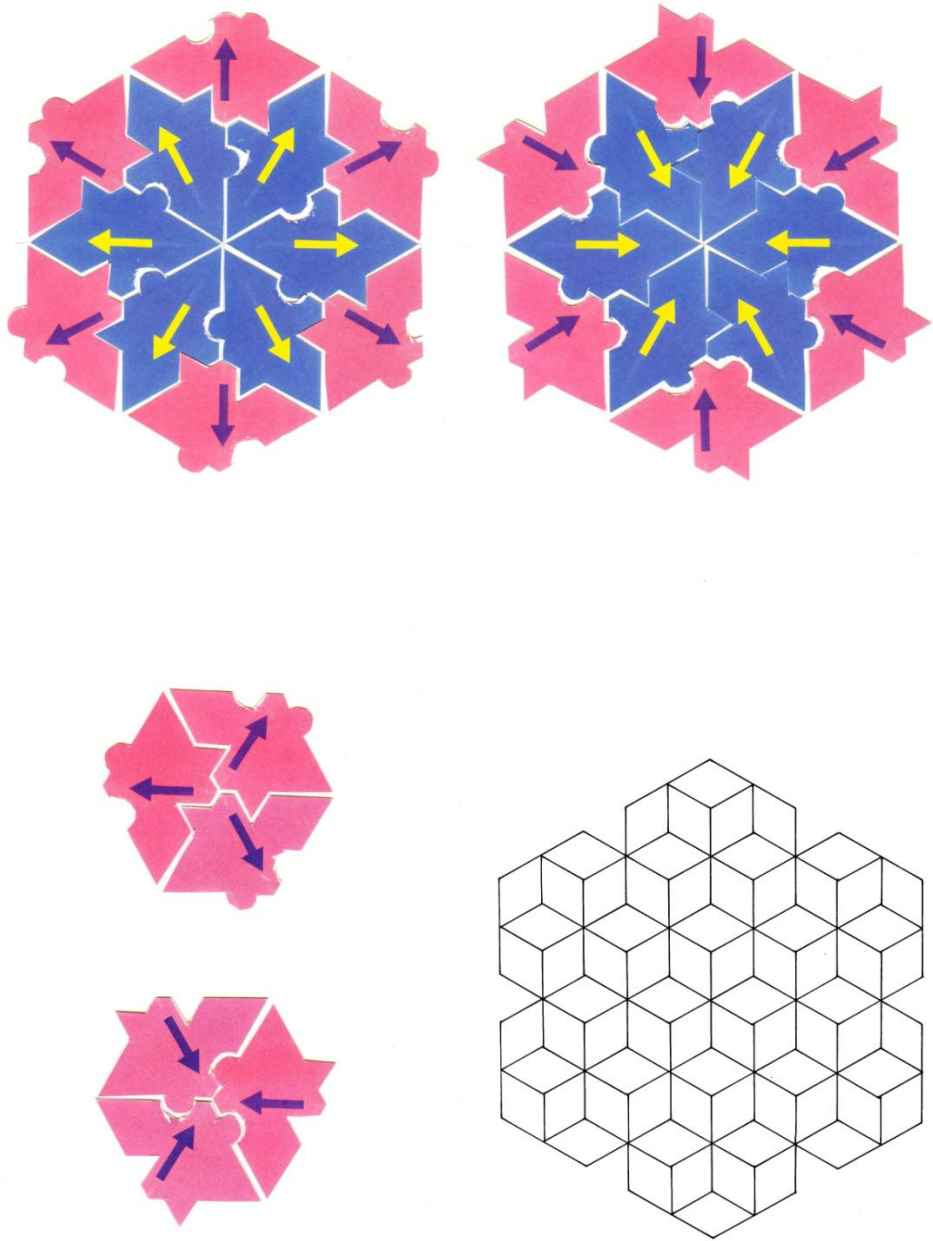

Fig. 2: The notches of the acute diamonds fit into two types of 6-fold stars. The obtuse diamonds are needed to complete a hexagon structures. The acute diamonds fit into two types of hexagons with 3-fold rotational symmetry. The underlying lattice is a periodic hexagonal tessellation.

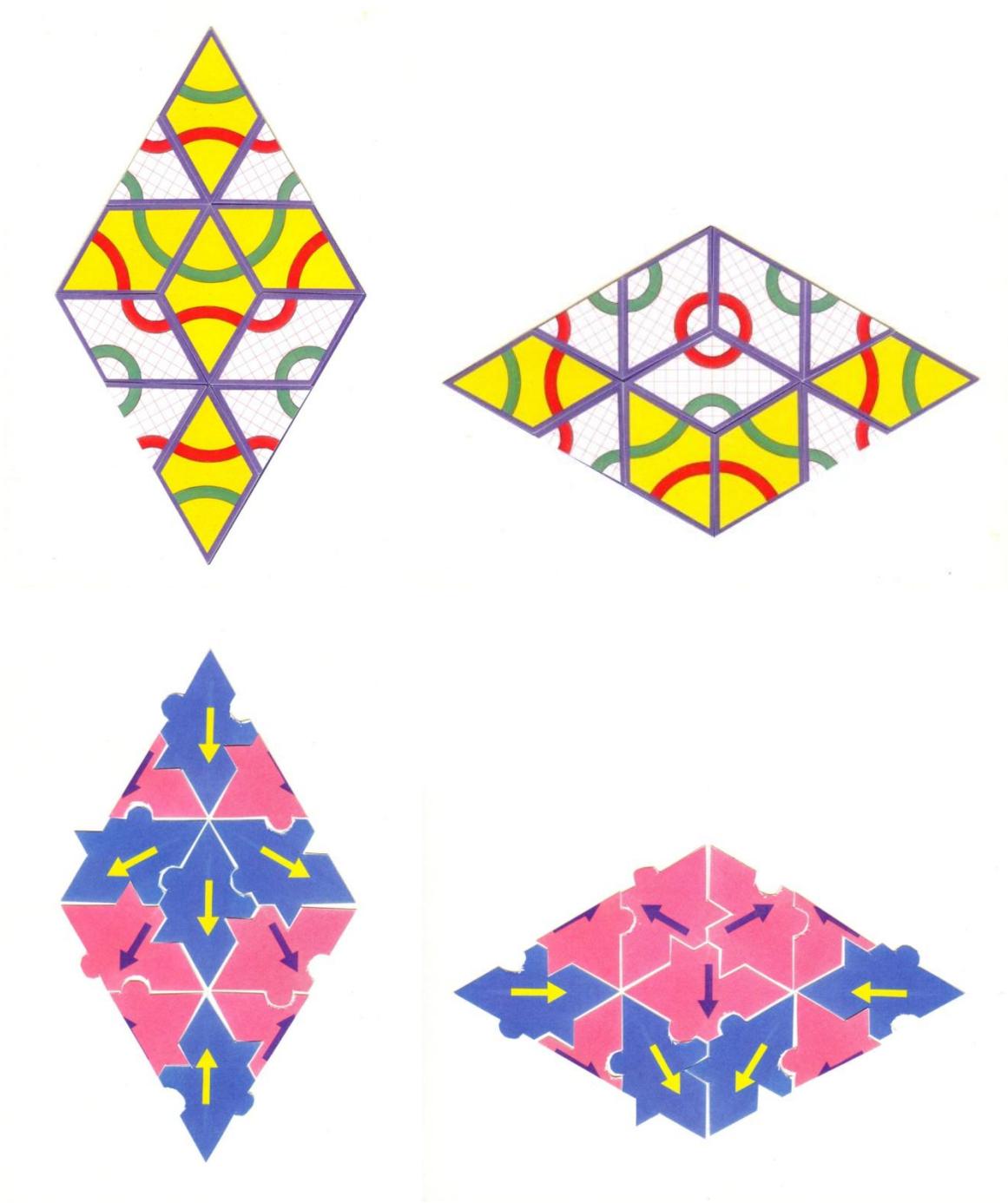

Fig. 3: An inflation scheme shown for the two embodiments of the two base diamonds: the inflated acute and obtuse diamonds then serve as the new base shapes for the next generation.

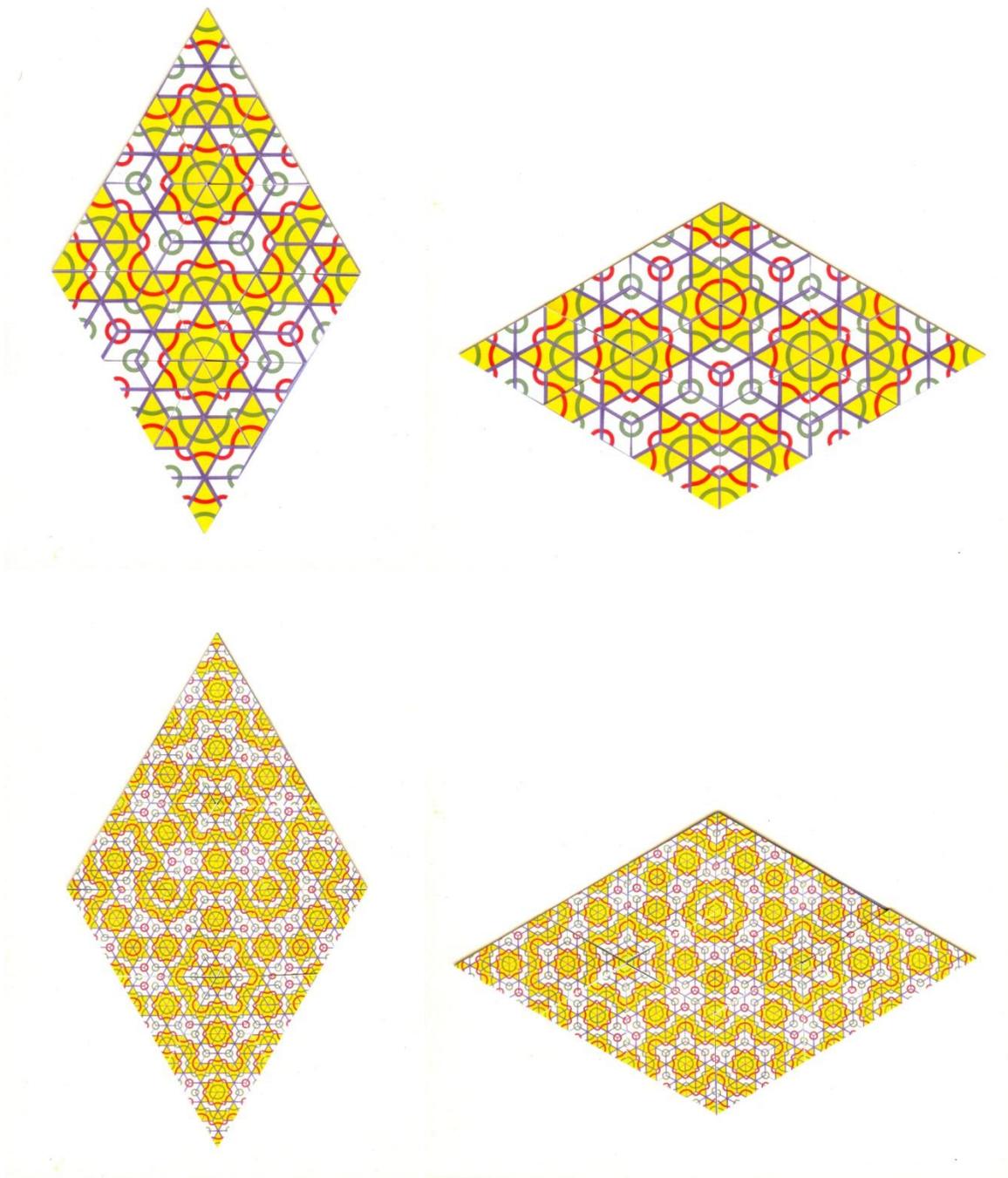

Fig. 4: Generation 3 (above) and 4 (below) of the inflation scheme using 2-colored bars, leading to new acute and obtuse base shapes that can be inflated to higher generations to form an infinite quasi-periodic 6-fold tiling.

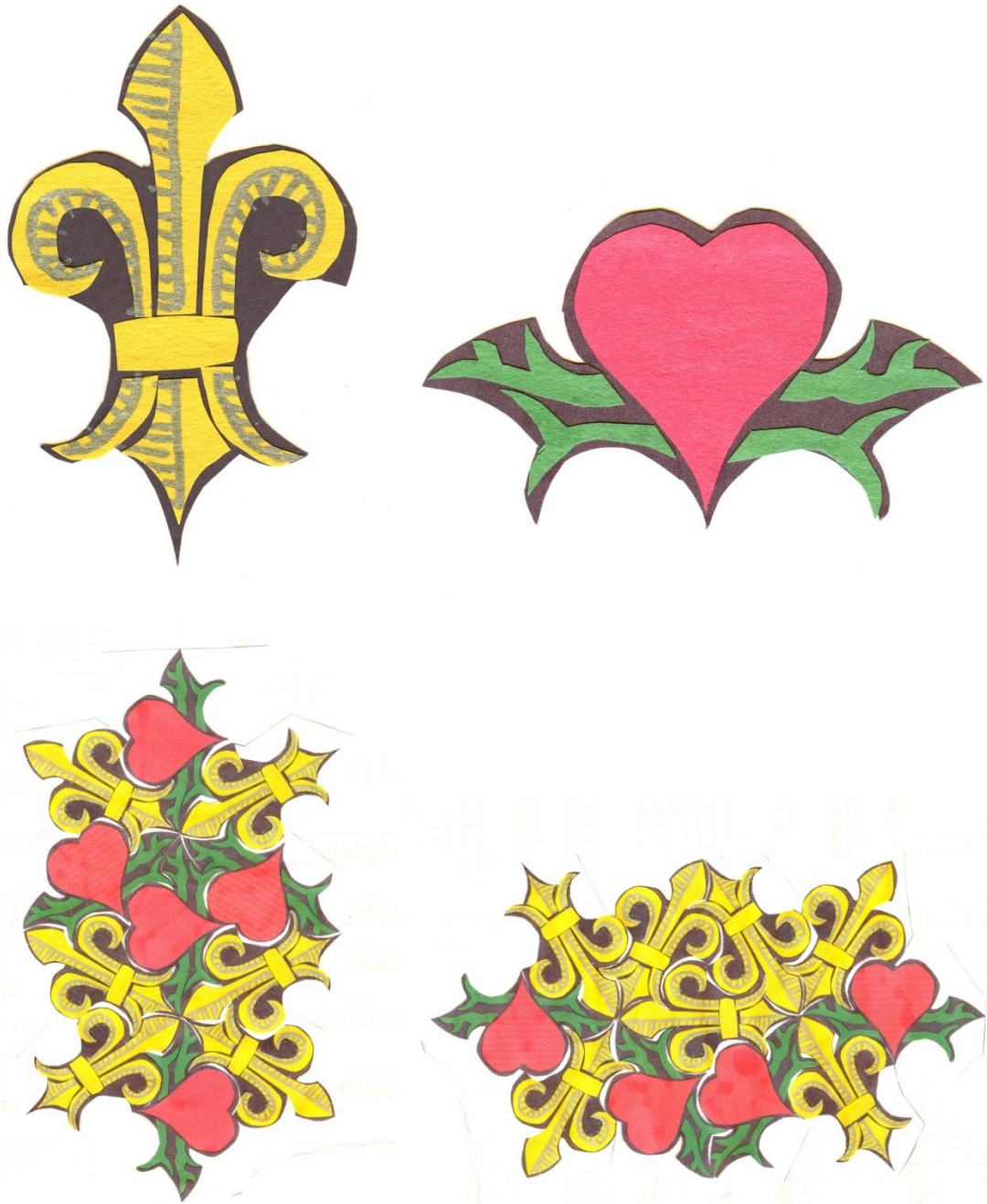

Fig. 5: Escher-esque figures of a Fleur-de-lis (deformed obtuse diamond) and a heart motif "Coeur d'Alene" (deformed acute diamond), shown rotated by 90 degrees (as compared to previous figures). At bottom generation 2 of the "fleurs" and "coeurs".

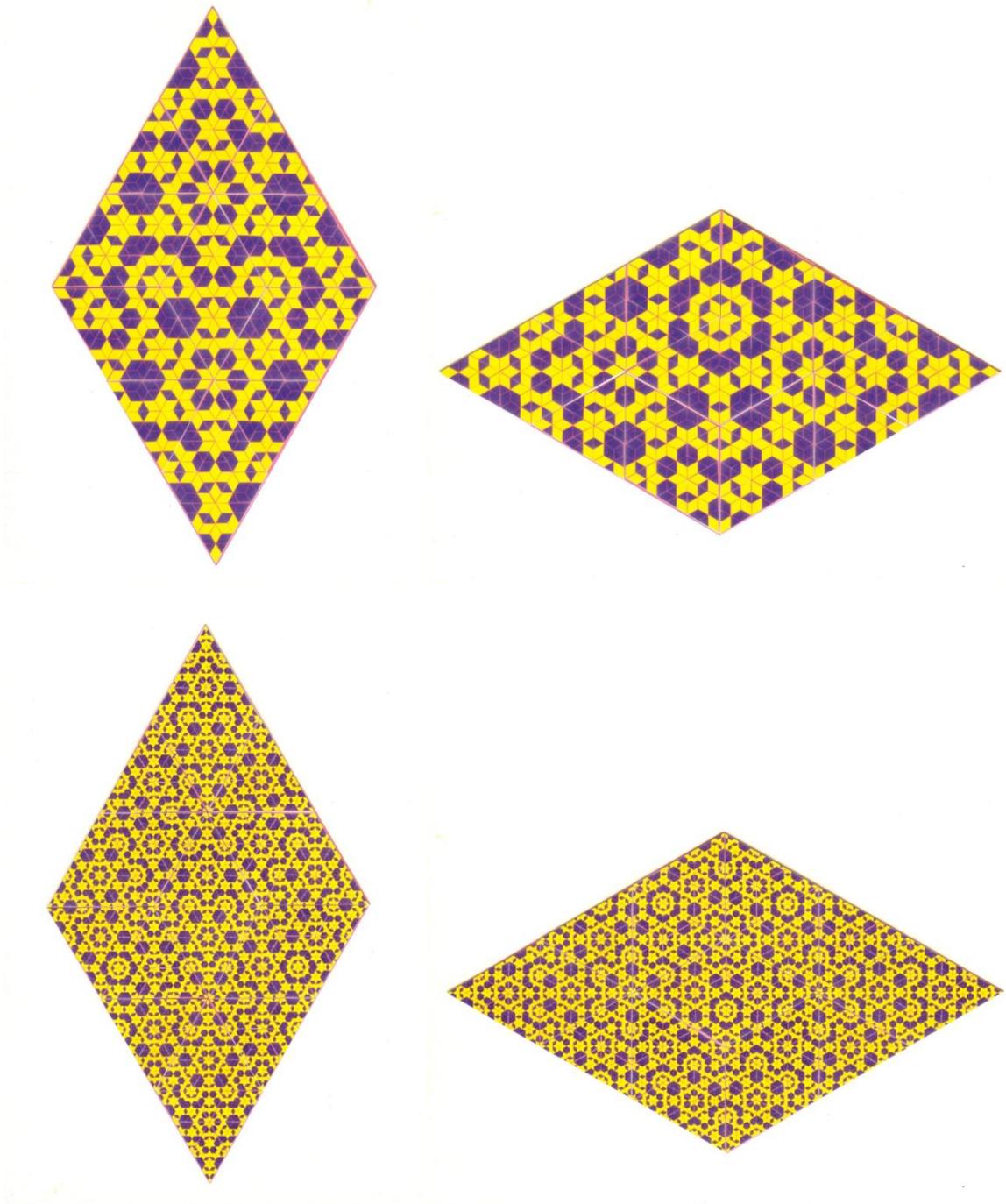

Fig. 6: Another embodiment of acute (gold) and obtuse diamonds (purple) with all markings (arrows, etc.) removed. Generation 4 and 5 of the inflation scheme.

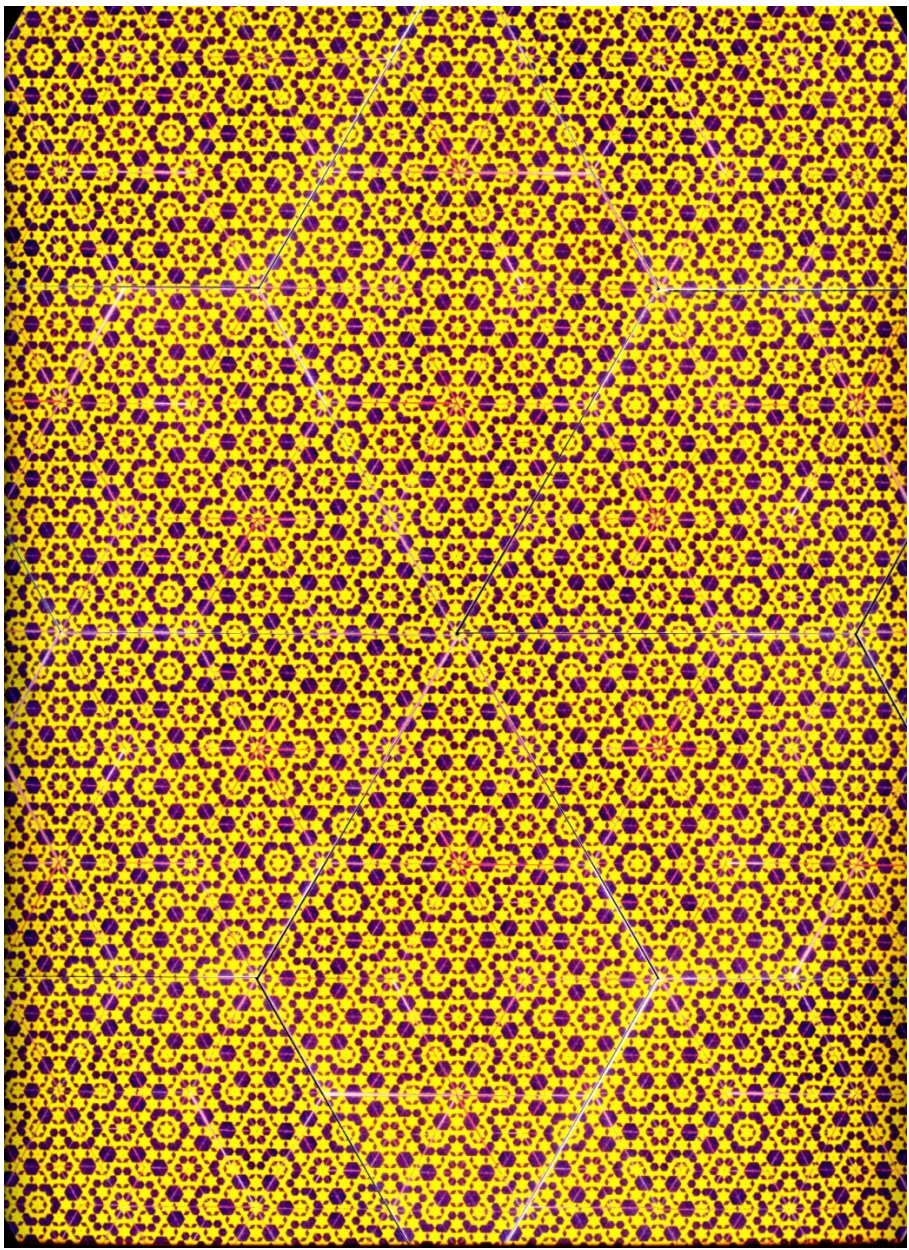

Fig. 7: Gold and purple diamonds in a 6-fold quasi-periodic tiling. Note the aperiodic structures spanning the entire scale of this embodiment.

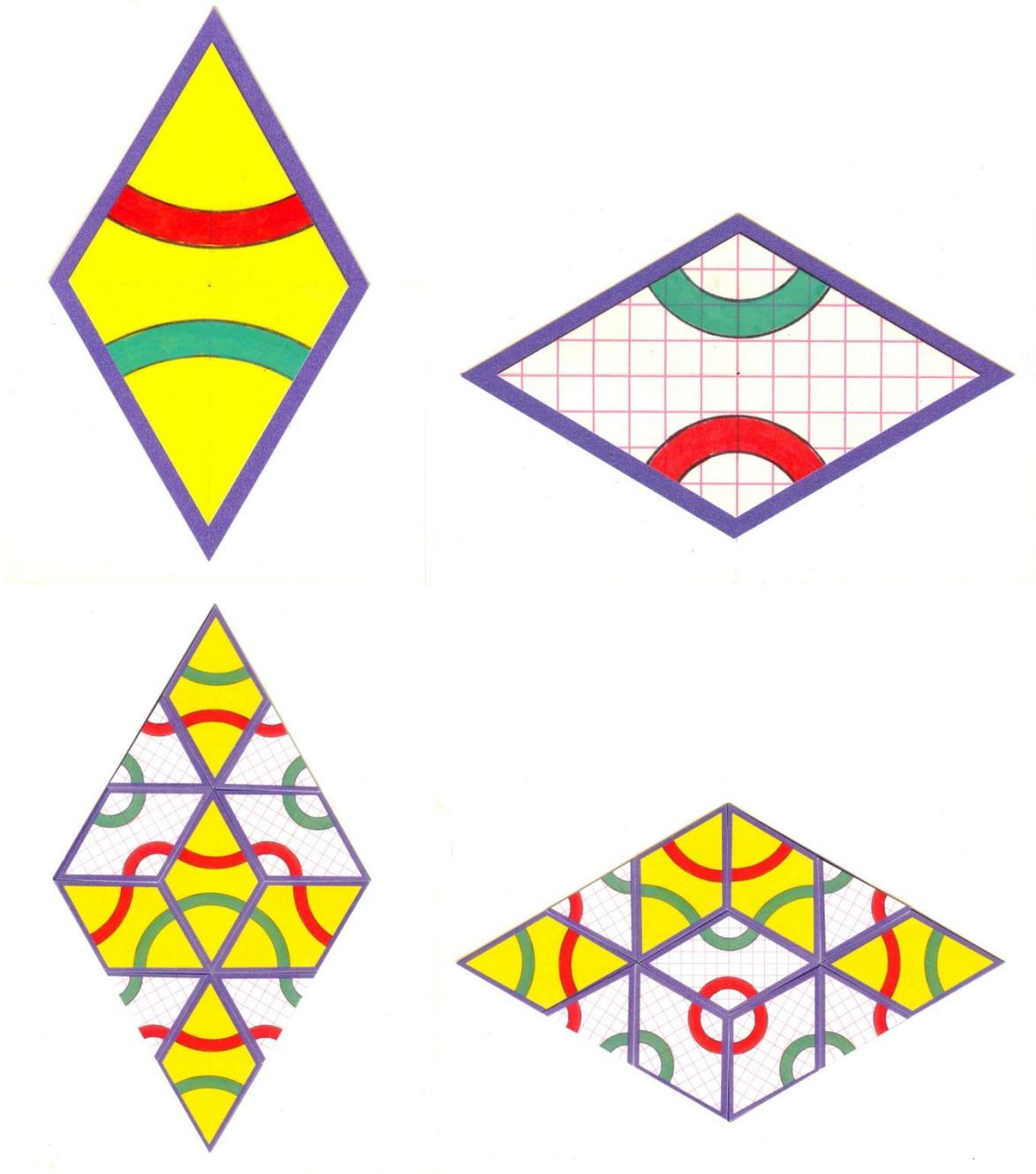

Fig. 8: A slightly different inflation scheme (variation 2) leading to a similar but slightly different 6-fold quasiperiodic tiling.

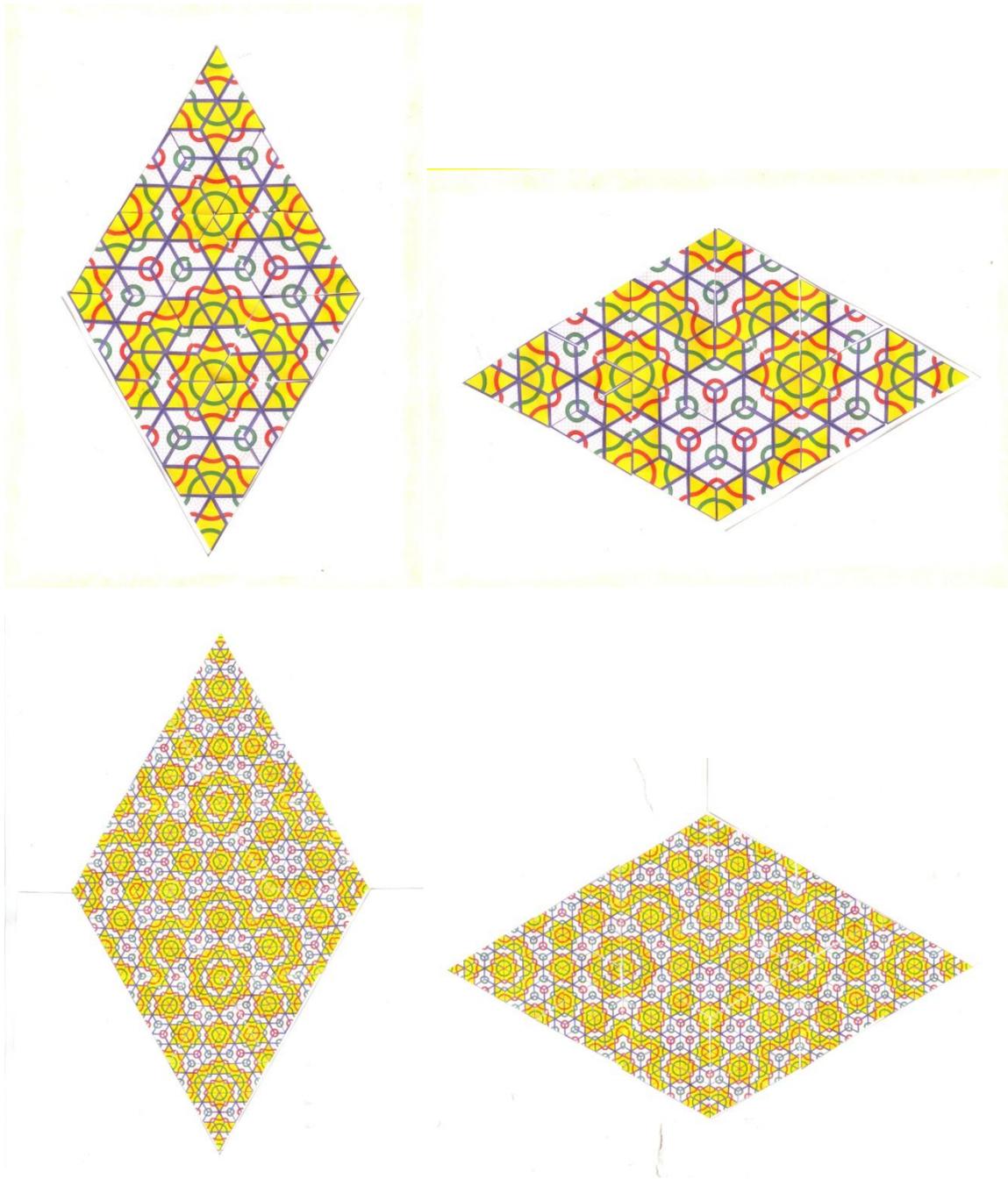

Fig. 9: Generation 3 and 4 of variation 2 of the inflation scheme.

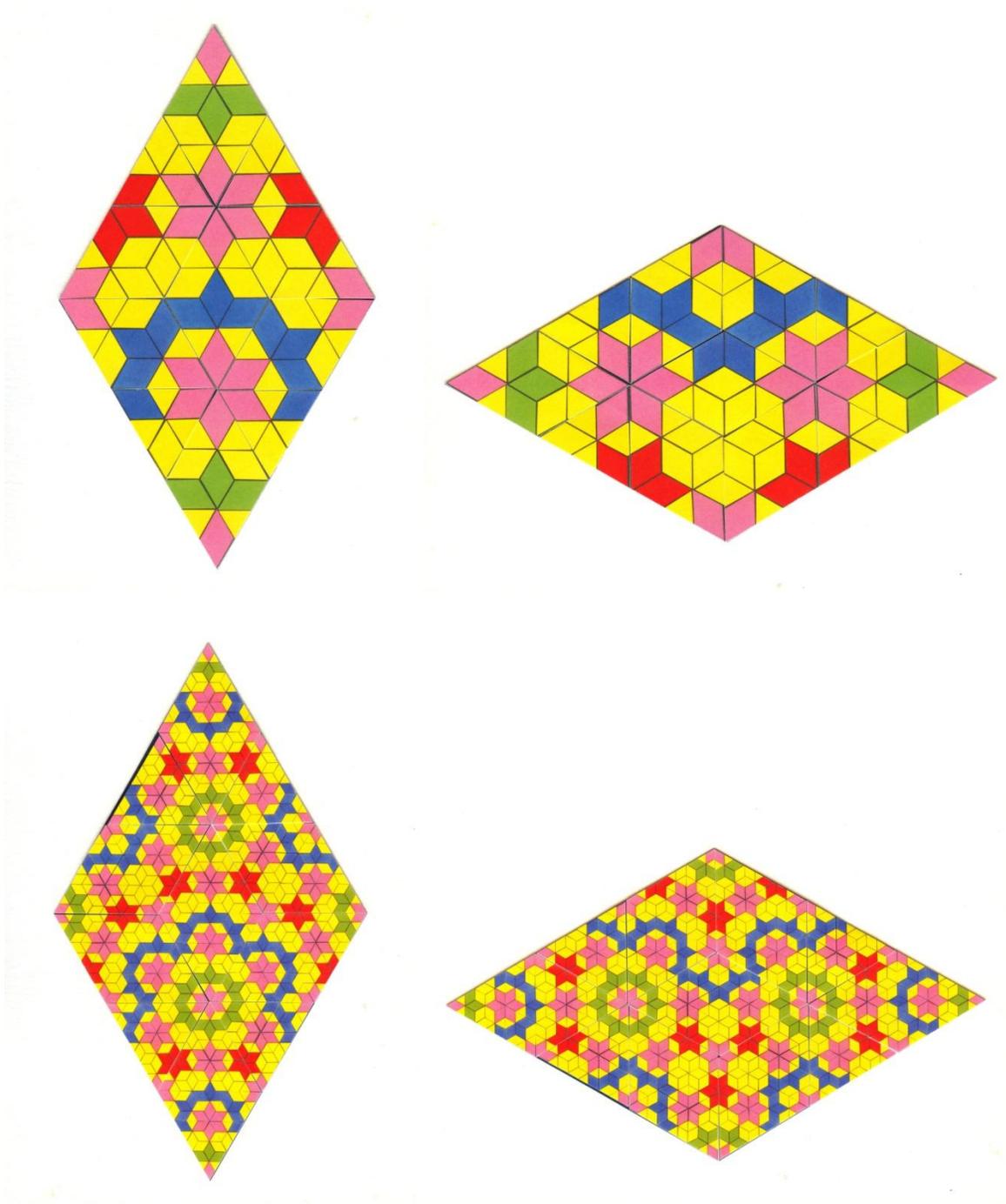

Fig. 10: Acute and obtuse diamonds, 5-colored, variation 2, generations 3 and 4. Note the periodic tessellation of pink stars and the periodic green rings in generation 4; the longest blue ribbon is aperiodic and ever longer in each new generation.

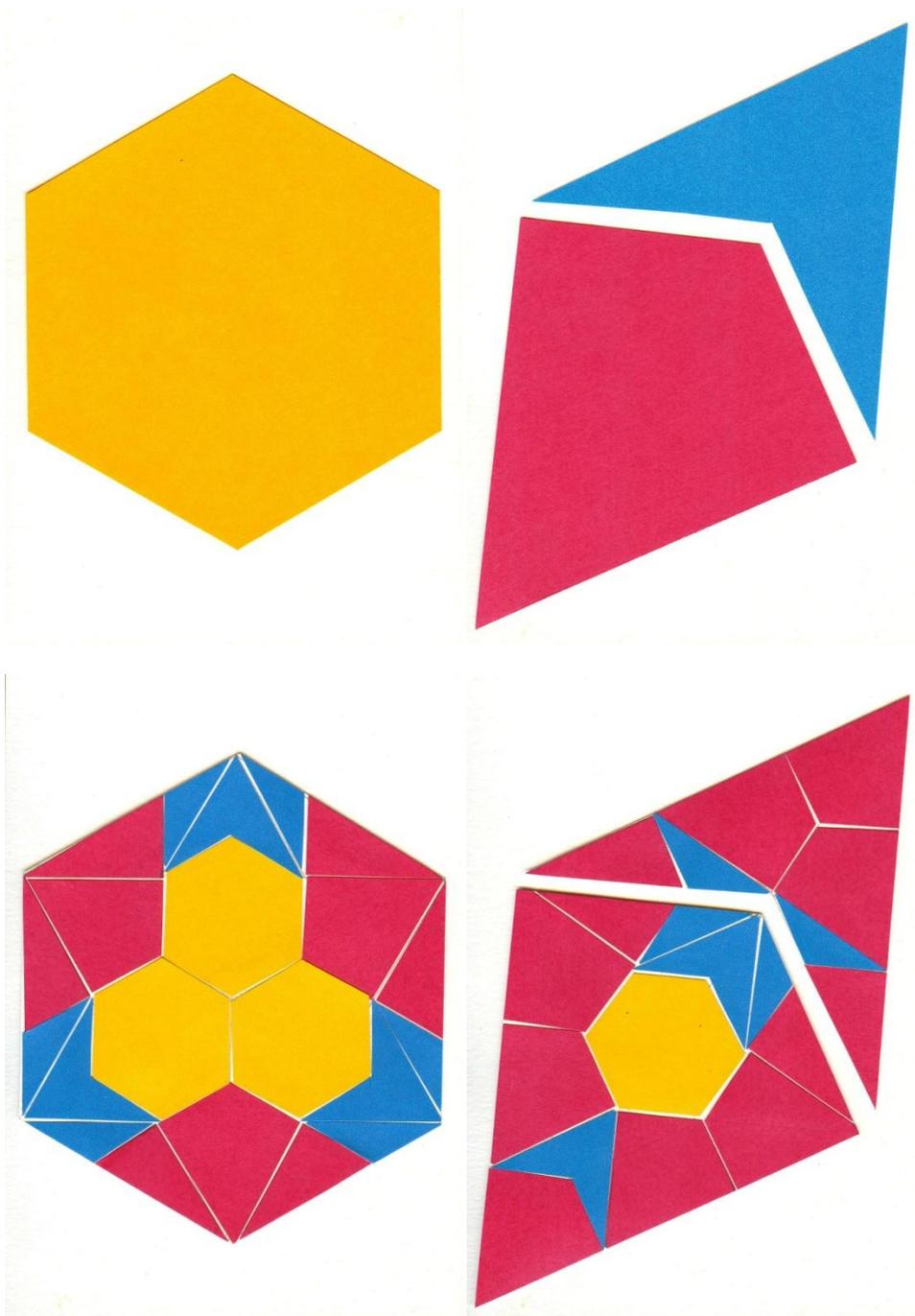

Fig. 11: Shield, kite, and dart. Substitution rule for the next generation. Again, there is also a slightly different substitution rule possible (variation 2, not shown).

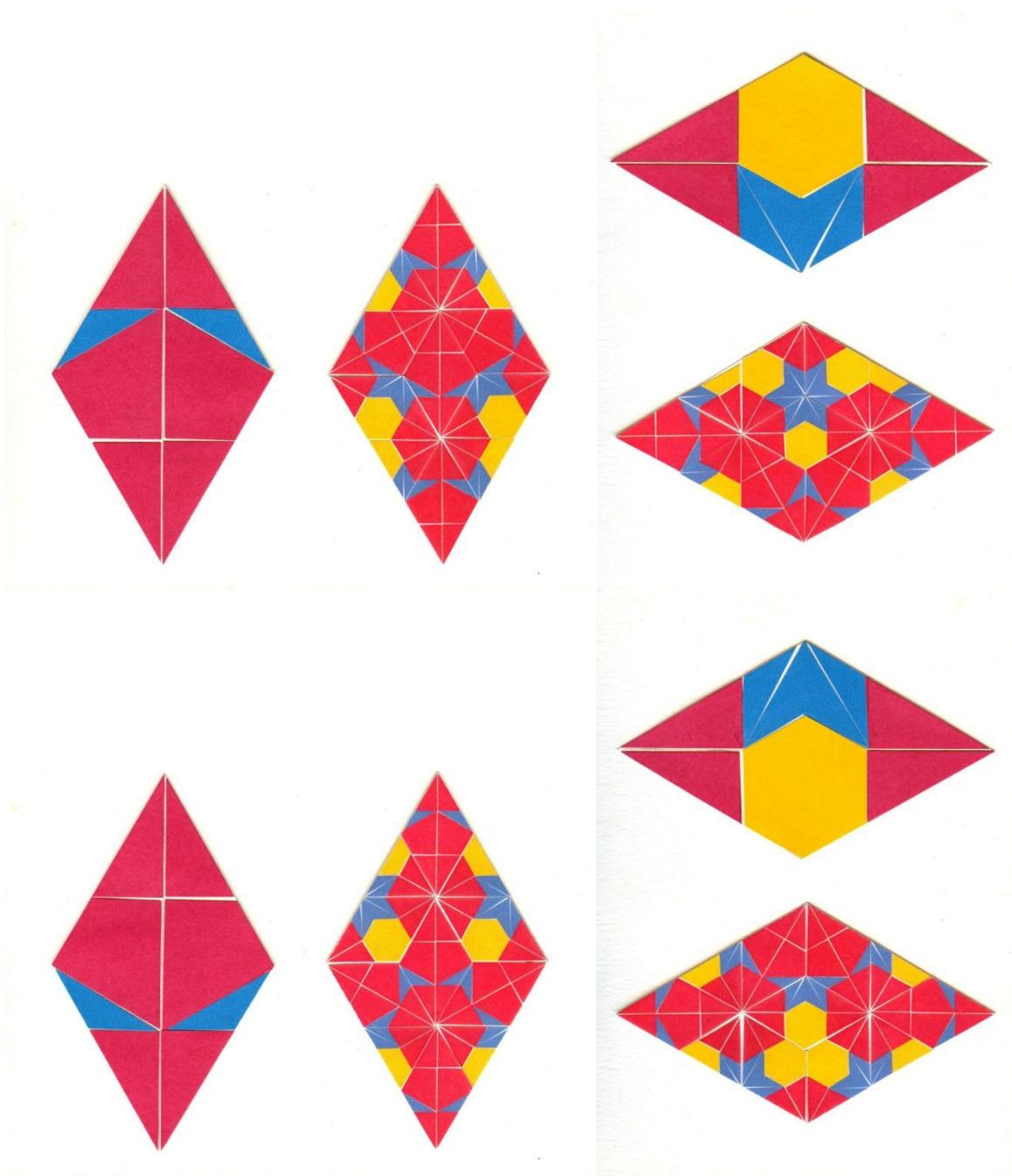

Fig. 12: Acute and obtuse diamonds tiled with darts, kites, and shields. Upper diamonds are variation 1, lower variation 2. Half-darts, half-kites, and half-shields must be matched with adjacent tiles to complete basic shapes.

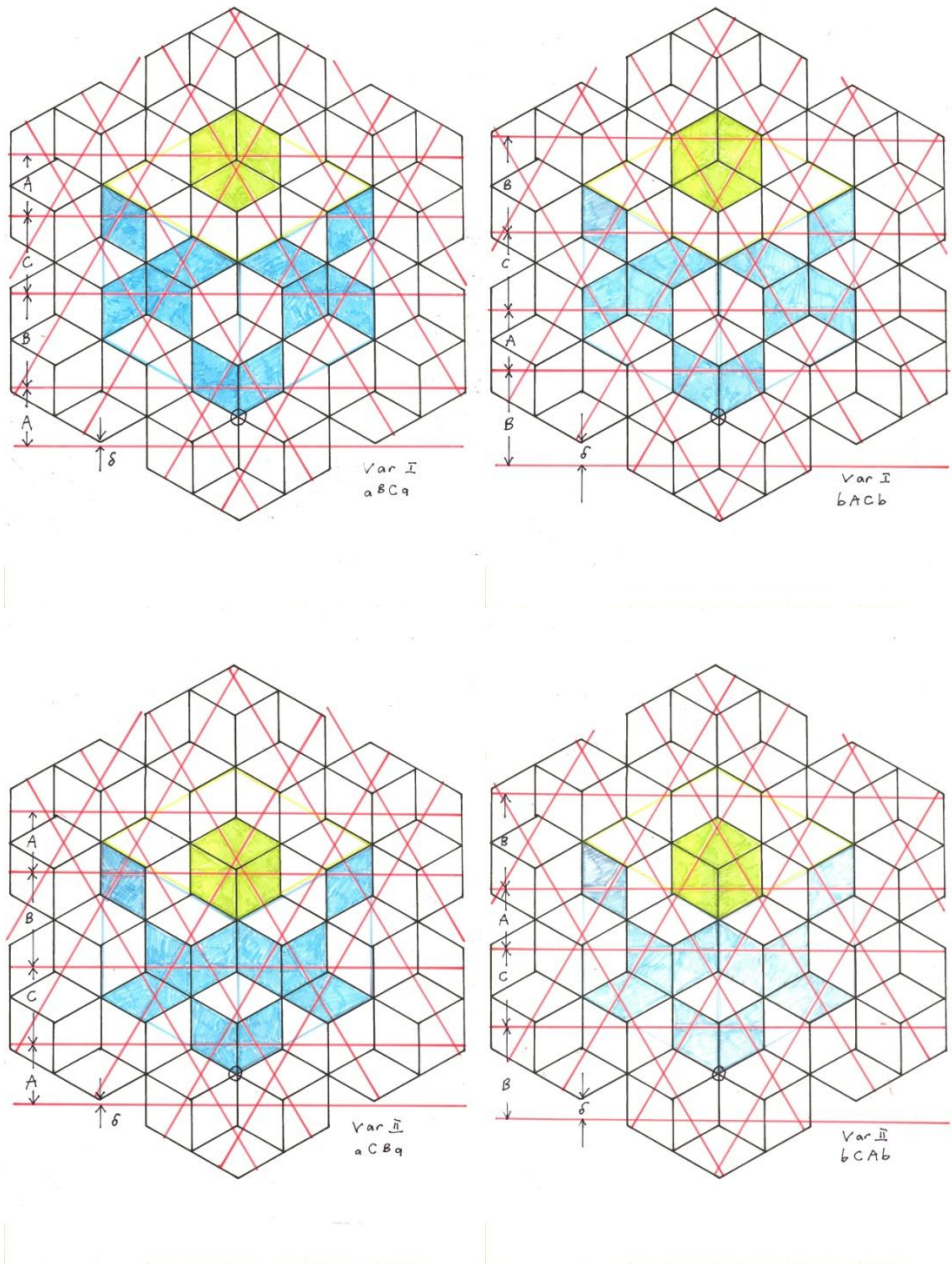

Fig. 13: Multigrid spacing for the two variations of the 6-fold quasiperiodic tilings. The circle marks the 6-fold center (see Fig. 2) where acute diamonds meet in a 6-fold star. Along the long diagonal of the inflated acute diamond, the pattern for variation 1 (2) is XYY (XXY), where X is the inverse of Y (see Fig. 3 & 8 for color scheme). Switching A with B simply provides a different polarity of the basic tiles without changing the overall substitution rule. The multigrid itself is an aperiodic pattern, dual to the diamond tiling. Also shown is the green shield of the obtuse diamond in the second generation.